\documentclass[10pt]{amsart}
\usepackage{amsthm}
\usepackage{amsmath}
\usepackage{amssymb}

\newcommand{\T}{\mathbb T}
\newcommand{\N}{\mathbb N}
\newcommand{\bee}{\begin{equation}}
\newcommand{\eee}{\end{equation}}

\newcommand{\Lb}{\mbox {\boldmath ${\Lambda}$}}

\def\ov{\overline}

\def\b0{{\bf 0}}

\newcommand{\be}{\begin{eqnarray}}
\newcommand{\ee}{\end{eqnarray}}
\newcommand{\supp}{\mbox{\rm supp}}

\newcommand{\R}{{\mathbb R}}

\newcommand{\Z}{{\mathbb Z}}

\newcommand{\Ak}{{\mathcal A}}

\newcommand{\Dk}{{\mathcal D}}

\newcommand{\Tk}{{\mathcal T}}

\newcommand{\Lam}{{\Lambda}}

\newtheorem{theorem}{Theorem}[section]
\newtheorem{lemma}[theorem]{Lemma}

\newtheorem{cor}[theorem]{Corollary}

\theoremstyle{definition}

\newtheorem{remark}[theorem]{Remark}

\numberwithin{equation}{section}

\begin{document}

\title{Strong coincidence and Overlap coincidence}

\address{Institute of Mathematics, University of Tsukuba\\
Tennodai 1-1-1, Tsukuba, Ibaraki\\
305-8571 Japan}
\email{akiyama@math.tsukuba.ac.jp}
\author{Shigeki Akiyama}
\date{}
\thanks{
This research was supported by the
Japanese Society for the Promotion of Science (JSPS), Grant in aid
21540012.}
\maketitle

\begin{abstract}
We show that strong coincidences of a certain many choices of
control points are equivalent to overlap coincidence for
the suspension tiling of Pisot substitution.
The result is valid for degree $\ge 2$ as well, under certain topological
conditions.
This result
gives a converse of the paper \cite{Akiyama-Lee:13}
and elucidates the tight relationship between two coincidences.
\end{abstract}

\section{Introduction}

Self-affine tiling dynamical system in $\R^d$
is a generalization of {\it substitution} dynamical system
on letters, which gives a nice model of self-inducing structures appear
in dynamical systems, number theory and the mathematics of
aperiodic order. Pure discreteness of self-affine tiling dynamics is
long studied from many points of views.
The idea of {\it coincidence}\footnote{In their notation, the {\it column number} is one.} appeared firstly in Kamae \cite{Kamae:72}, and then in a comprehensive form in Dekking \cite{Dekking:78} for constant length substitution (see also \cite{Queffelec:87}).
Generalizing a pioneer work of Rauzy \cite{Rauzy:82},
Arnoux-Ito \cite{Arnoux-Ito:01} gave
a geometric realization of irreducible Pisot unit substitution of degree $d$.
They defined {\it strong coincidence},
which ensures that
their geometric substitution gives rise to a domain exchange of $\R^{d-1}$,
which is also semi-conjugate to the toral rotation of $\T^d$.
It is remarkable that in many cases, it is even conjugate to the total rotation,
which immediately implies that the system in pure discrete (see \cite{
Barge-Diamond:02,Ito-Rao:03,
Barge-Kwapisz:06, Siegel-Thuswaldner:09} for further developments).
On the other hand, {\it overlap coincidence}
introduced by Solomyak \cite{Solomyak:97}
is an equivalent condition for pure discreteness of a given self-affine
tiling dynamical system.
This is also described
as a geometric/combinatorial condition which guarantees that the tiling
and its translation by return vectors become exponentially close if we
iteratively enlarge return vectors by substitution. Lee \cite{Lee:07}
showed deeper characterizations
that overlap coincidence is equivalent to {\it algebraic coincidence},
and the fact that the corresponding point set is an inter-model set.

Until now, the relation between
strong coincidence and overlap coincidence is not fully understood.
Motivated by the claim of Nakaishi \cite{Nakaishi:11}, Akiyama-Lee
\cite{Akiyama-Lee:13} generalized the notion of strong coincidence
to $\R^d$ and
showed that overlap coincidence implies strong coincidence,
and moreover {\it simultaneous coincidence},
provided that the associated point set is admissible and its height group is trivial.
In this paper we shall give a converse statement
for the suspension tiling of Pisot substitution
at the expense of
assuming many strong coincidences at a time, that is,
strong coincidences on a certain many choices of
control points imply overlap coincidence and vice versa.
If every tile is connected and the tiling is not a collection of
unbounded connected identical colored patches,
then the same result holds for $d\ge 2$ (Theorem \ref{Main}).
This result elucidates the tight
relationship between two coincidences.

\section{Terminologies}
\subsection{Tiles and tilings}

We shall briefly recall basic definitions used in this paper.
A {\em tile} in $\R^d$ is defined as a pair $T=(A,i)$ where $A$
is a compact
set in $\R^d$ which is the closure of its interior, and
$i=\ell(T)\in \{1,\ldots,m\}$
is the {\it color} of $T$. We call $A$ the {\it support} of $T$ and denote $\supp(T)=A$.
The {\it translate} of $T$ is defined by $g+T = (g+A,i)$ for $g \in \R^d$.
Let $\Ak = \{T_1,\ldots,T_m\}$ be a finite set of tiles in
$\R^d$ such that $T_i=(A_i,i)$; we will call them {\em
prototiles}. A tiling $\Tk$ is a collection of translates of prototiles
which covers $\R^d$ without interior overlaps.
A finite collection of tiles which appear in $\Tk$ is called a {\em patch}.
A {\it generalized patch} is a collection of tiles in $\Tk$ whose
cardinality is not necessarily finite.
Its support is defined to be the union of the supports of tiles.
The diameter of a generalized patch is the supremum of Euclidean distance of two
points lie within the support of the patch.
A map $\Omega$ from $\Ak$ to the set of patches is called
a {\em substitution} with a
$d\times d$ expansive matrix $Q$ if there exist finite sets $\Dk_{ij}\subset
\R^d$ for $i,j \le m$ such that
\begin{equation}
\Omega(T_j)=
\{T_i+u:\ u\in \Dk_{ij},\ i=1,\ldots,m\}
\label{subdiv}
\end{equation}
with
\begin{eqnarray} \label{tile-subdiv}
Q A_j =
\bigcup_{i=1}^m (A_i+\Dk_{ij})
=\bigcup_{i=1}^m \bigcup_{u\in \Dk_{ij}} (A_i+u) \ \ \ \mbox{for} \ j\le m,
\end{eqnarray}
and the last union has mutually disjoint interiors.
The substitution (\ref{subdiv}) extends to all
translates of prototiles and patches in a natural way.
A {\it substitution tiling} of $\Omega$ is a tiling $\Tk$
that all the patches of $\Tk$ is a sub-patch of
$\Omega^n(T)$ for some $n\in \N$ and $T\in \Tk$.
A substitution tiling $\Tk$ is a {\em fixed point} of $\Omega$ if $\Omega(\Tk) = \Tk$
holds. We say that a substitution tiling is {\em primitive} if
the corresponding substitution matrix $M=(\sharp \Dk_{ij})$ is primitive, and
{\em irreducible} if the characteristic polynomial of
$M$ is irreducible.
We say that $\Tk$ has {\em finite local complexity} (FLC)
if for any $R$ there are only finitely many patches
of diameter less than $R$ up to translation.
A tiling $\Tk$ is {\em repetitive} if every
patch is relatively dense in $\Tk$.
A FLC substitution tiling
of a primitive substitution is called a {\em self-affine tiling}.
Every self-affine tiling is repetitive, which follows from the primitivity of
substitution.
Let $\lambda>1$ be the Perron-Frobenius eigenvalue  of the substitution matrix $M$
and $D$ be the set of eigenvalues of $Q$. By the tiling criterion of
Lagarias-Wang \cite{Lagarias-Wang:03}, $\lambda$ is the element of $D$
of maximum modulus.
We
say that $Q$ fulfills {\em Pisot family condition} if every algebraic conjugate
$\mu$ of an element of $D$ with $|\mu|\ge 1$ is contained in $D$.
\medskip

The set of all substitution tilings of $\Omega$ forms a {\it tiling space}.
By using a fixed point $\Tk$ of $\Omega$, we can describe this space
as the orbit closure of $\Tk$ under the translation action: $X_{\Tk} =
\ov{\{\Tk-g :\,g \in \R^d \}}$, the closure is taken by `local topology'.
The FLC assumption implies
$X_{\Tk}$ is compact and we get a topological dynamical system
$(X_{\Tk},\R^d)$ where $\R^d$ acts by translations.
This system is
minimal and uniquely ergodic (\cite{Solomyak:97, LMS}), and
we are interested in the spectra of self-affine tiling dynamical systems.
Tiling dynamical system $X_{\Tk}$
has {\em pure discrete spectrum} if the eigenfunctions for the $\R^d$-action
forms a complete orthonormal basis of $L^2(X_{\Tk},\mu)$ \cite{Solomyak:97}.

\subsection{Control points}

\noindent
A Delone set is a
relatively dense and uniformly discrete subset of $\R^d$. We say
that $\Lb=(\Lambda_i)_{i\le m}$ is a {\em Delone multi-color set}
in $\R^d$ if each $\Lambda_i$ is Delone and
$\cup_{i=1}^m \Lambda_i \subset \R^d$ is Delone.
We say that $\Lam \subset \R^d$ is a {\em Meyer set}
if it is a Delone set and $\Lam - \Lam$ is uniformly discrete in $\R^d$
\cite{Lagarias:96}.
$\Lb = (\Lam_i)_{i\le m}$ is called a {\em
substitution Delone multi-color set} if $\Lb$ is a Delone multi-color set and
there exist an expansive matrix
$Q$ and finite sets $\Dk_{ij}$ for $i,j\le m$ such that
\begin{equation}
\label{eq-sub}
\Lambda_i = \bigcup_{j=1}^m (Q \Lambda_j + \Dk_{ij}),\ \ \ i \le m,
\end{equation}
where the union on the right side is disjoint.

Given a fixed point $\Tk$ of $\Omega$,
we can associate a substitution Delone multi-color set
$\Lb_{\Tk} = (\Lam_i)_{i \le m}$ of $\Tk$ by
taking representative points of tiles in the relatively same positions for the same color tiles in the tiling.
There is a canonical way to choose representative points, called {\it control points}.
A tile map $\gamma=\gamma_{\Omega}$ is a map from $\Tk$ to itself
which sends a tile $T$ to the one in $\Omega(T)$ such that $\gamma(T_1)$ and $\gamma(T_2)$ are located in the same relative position
in $\Omega(T_1)$ and $\Omega(T_2)$ whenever $\ell(T_1)=\ell(T_2)$.
A control point $c(T)$ of $T\in \Tk$ is defined by
$$
c(T)=\bigcap_{n=1}^{\infty} Q^{-n} (\gamma^n(T)).
$$
Control points are representative points, i.e., $U-c(U)=V-c(V)$ holds
if $\ell(U)=\ell(V)$ with $U,V\in \Tk$.
Let $\Lam_i$ be the set of control points of color $i$. Clearly, $j=\ell(\gamma(T_i))$ implies $Q \Lam_i \subset \Lam_j$
and
the set of control points $\mathcal{C}=\bigcup_{i=1}^m \Lam_i$ is invariant under the
expansion by $Q$, that is, $Q \mathcal{C} \subset \mathcal{C}$.
We obtain an associated substitution Delone multi-color set
$\Lb = \Lb_{\Tk} = (\Lam_i)_{i \le m}$.

In section \ref{Coincidences}, we have to assume
a lot of strong coincidences
by changing control points for a given tiling $\Tk$.
If we change control points of tiles of $\Tk$
by $\Lam'_i=\Lam_i-g_i$, then the set equation will be shifted like
$$
\Lam'_i= \bigcup_{j=1}^m Q \Lam'_j + \mathcal{D}'_{ij}
$$
with $\mathcal{D}'_{ij}=\left\{ d_{ij}+Q g_j-g_i\ :\ 1 \le i,j \le m \right\}$.
The corresponding tile equation becomes
$$
Q A'_j= \bigcup A'_i + \mathcal{D}'_{ij}
$$
which is satisfied by $A'_j=A_j+g_j$. So we set $\supp(T'_j)=A'_j$ and
$\ell(T'_j)=\ell(T_j)$.
To avoid heavy notation, we do not
distinguish such changes of control points and use the
same symbols $\Lam_i$ and $T_i$.


\subsection{Coincidences}
The set of {\it return vectors} is defined by
$\Xi(\Tk) = \{y \in \R^d \ : \ U = V-y, \ \mbox{where $U, V \in \Tk$}\}$.
A triple $(U, y, V)$, with $U, V \in \Tk$
and $y \in \Xi(\Tk)$, is called an {\em overlap} if
\[ (\supp(U))^{\circ} \cap (\supp(V)-y)^{\circ} \neq \emptyset. \]
An overlap $(U, y, V)$ is a {\em coincidence} if $U = V-y$.
Let $\mathcal{O} = (U, y, V)$ be an overlap in $\Tk$, we define
{\em $\ell$-th inflated overlap}
\begin{eqnarray*}
{\Omega}^{\ell} \mathcal{O} = \{(U', Q^{\ell} y, V') \, :
U' \in \Omega^{\ell}(U), V' \in \Omega^{\ell}(V), \ \mbox{and $(U',Q^{\ell}y,V')$
is an overlap} \}.
\end{eqnarray*}
We say that a self-affine tiling $\Tk$ admits {\em overlap
coincidence} if there exists $\ell \in \Z_+$ such that for each
overlap $\mathcal{O}$ in $\Tk$, ${\Omega}^{\ell} \mathcal{O}$
contains a coincidence.
Two overlaps $(U,y,V)$ and $(U_1,y_1,V_1)$ are equivalent, if
there is $x\in \R^d$ that both $U_1=U-x$ and $V_1-y_1=V-y-x$ hold.
The equivalence class is denoted by $\widetilde{(U,y,V)}$.
Hereafter we assume an important condition that
$\Xi(\Tk)$ forms a Meyer set. This condition
is equivalent to the Pisot family condition
for $Q$, if $Q$ is diagonalizable and all its eigenvalues are
algebraic conjugate with the same multiplicity \cite{Lee-Solomyak:12}. The number of equivalence classes of
overlaps is finite, by the Meyer property of $\Xi(\Tk)$.
The action of $\Omega$ is well-defined on equivalence classes of overlaps.
An {\it overlap graph with multiplicity} is
a finite directed graph whose vertices are
the equivalence classes of overlaps. Multiplicities of the edge from
$\widetilde{(U,y,V)}$ to $\widetilde{(A,z,B)}$ is
given by the number of overlaps in $\Omega((U,y,V))$
equivalent to $(A,z,B)$ (c.f. \cite{Akiyama-Lee:10}).
Overlap coincidence is confirmed by checking whether from each vertex
of this graph there is a path leading to a coincidence.
Overlap coincidence is equivalent to
pure discreteness of self-affine tiling dynamical system $X_{\Tk}$
\cite{Solomyak:97}.

Strong coincidence on letter substitution is naturally generalized to
self-similar tiling in $\R^d$ in \cite{Akiyama-Lee:13}.
We adapt this definition to control points.
Let $\Tk$ be a self-affine tiling in $\R^d$ and $\mathcal{A} = \{T_1, \cdots, T_m\}$ be the prototile set of $\Tk$.
We say that the set of the control points
is {\it admissible} if $\cap_{i \le m} (\supp(T_i)-c(T_i))$ has non-empty interior.

Let $\Tk$ be the fixed point of $\Omega$.
Let $c(T_i)\ (i=1,\dots,m) $ be the admissible control points and
$\Lb$ be an associated substitution Delone multi-color set
for which $\Tk = \{T_i -c(T_i)+u_i \ |  \ u_i \in \Lam_i, i \le m \}$.
If for any $1 \le i, j \le m$, there is a positive integer $L$ that
\begin{eqnarray} \label{StronCoin-I}
\Omega^L(T_i - c(T_i)) \cap
\Omega^L(T_j - c(T_j)) \neq \emptyset,
\end{eqnarray}
then we say that $\Lb$ admits {\em strong coincidence}.
In other words, strong coincidence means that
for every pair of tiles $(U,V)\in \Tk^2$, $\Omega^L(U - c(U))$ and $\Omega^L(V - c(V))$
share a common tile in the same position for some $L$.

\section{Strong coincidence and Overlap coincidence}
\label{Coincidences}

The set of {\it eventually return vectors} is defined by
$$
\mathcal{G} :=\bigcup_{k=0}^{\infty} Q^{-k} (\Lam_i-\Lam_i), \ \ \  \mbox{for some $i \le m$}
$$
which is independent of the choice of $i$, by primitivity of $\Omega$.
The tiling dynamical system is invariant under replacement of the substitution rule
$\Omega$ by $\Omega^n$. We consider control points of $\Omega^n$ as well.
Hereafter we put $\Lam=\bigcup_{i=1}^m \Lam_i$ for
$\Lb=\Lb_{\Tk}=(\Lam_i)$ to distinguish the multi-color set and its union.
Let $\langle \mathcal{G} \rangle$ be the additive subgroup of $\R^d$
generated by $\mathcal{G}$. We say that
$\Tk$ satisfies {\it multiple
strong coincidence of level $n$} if all multi-color Delone set
$\Lb$'s generated by admissible control points of $\Omega^n$ with $\Lam -\Lam \subset
\langle \mathcal{G} \rangle$ admit strong coincidence.

Hereafter when
we speak about a topological/metrical property (connected, bounded, diameter)
of a generalized patch, it refers to the corresponding property of its support.
A {\it rod} is an unbounded connected generalized
patch of $\Tk$ whose tiles have an identical color.
A {\it rod tiling} is a tiling that every tile belongs to a rod.
For ease of negation,
a {\it non-rod } tiling is a tiling which is not a rod tiling.
A tiling is called {\it non-periodic} if there are no non-trivial period, i.e., $\{p\in \R^d\ |\ \Tk+p=\Tk\}=\{0\}$.

\begin{remark}
\label{RodTiling}
There are many examples of periodic self-affine rod tiling. Consider a tiling
of $\R^2$ by squares $[0,1]^2+(x,y)$ with $(x,y)\in \Z^2$ and their colors are
defined by $y \pmod 2$ or $x+y \pmod 2$. However
we do not know an example of non-periodic self-affine rod tiling.
\end{remark}

\begin{theorem}
\label{Main}
Let $\Tk$ be a non-rod self affine tiling by connected tiles
such
that $\Xi(\Tk)$ is a Meyer set.
Then there is a constant
$n$ depending only on $\Tk$ that
$\Tk$ satisfies multiple strong coincidence of level $n$ if and  only if
$\Tk$ satisfies overlap coincidence.
\end{theorem}

Consider a substitution $\sigma$ over $m$ letters $\{1,2,\dots,m\}$
whose substitution matrix is $M_{\sigma}=(|\sigma(j)|_i)$, where $|w|_i$
is the number of letter $i$ in a word $w$.
We say that $\sigma$ is a {\it Pisot} substitution, if the Perron
Frobenius root $\beta$ of $M_{\sigma}$ is a Pisot number.
The canonical {\em suspension tiling} $\Tk$ in $\R$ of $\sigma$
with an expansion factor $\beta$ is defined
by associating to the letters the intervals whose lengths are given by a left
eigenvector of $M_{\sigma}$ corresponding to $\beta$.

\begin{cor}
\label{1D}
The statement is valid for the suspension tiling of a Pisot substitution.
\end{cor}

Indeed, $1\times 1$ matrix $Q=(\beta)$ satisfies Pisot family condition,
tiles are intervals and the suspension tiling can not be a
rod tiling, since it has at least two translationally inequivalent tiles in $\R$.

\begin{remark}
\label{Small}
Multiple strong coincidence of level $n$ requires
many strong coincidences at a time for a fixed tiling $\Tk$ even when $n=1$.
In dimension one, the claim of Nakaishi \cite{Nakaishi:11} reads
a single strong coincidence implies overlap coincidence. Theorem
\ref{Main} covers general cases but the requirement is much stronger.
It would be interesting is to make smaller the constant
$n$ in Theorem \ref{Main}. For e.g., can we take $n=1$ ?
\end{remark}

We prepare a lemma.
\medskip

\begin{lemma}
\label{CycleExt}
Let $G$ be a strongly connected finite directed graph and $C$ be a set of
cycles of $G$. Then there is a subgraph $G(C)$ of $G$ with the following property.
\begin{itemize}
\item The set of vertices of $G(C)$ is equal to that of $G$.
\item Every vertex has exactly one outgoing edge.
\item The set of cycles of $G(C)$ is equal to $C$.
\end{itemize}
\end{lemma}

\proof
Put $H_0=C$.
We inductively
construct $H_i$ for $i=0,1,\dots$
which satisfies:
\begin{itemize}
\item Every vertex has exactly one outgoing edge.
\item The set of cycles of $H_i$ is equal to $C$.
\end{itemize}
Assume that the induced graph $G\setminus H_i$ is non empty and take
a vertex $v$ from $G\setminus H_{i}$.
Since $G$ is strongly connected, there is a path from $v$ leading to $H_{i}$.
So there is a vertex $u\in G\setminus H_i$ and an edge
from $u$ to a vertex of $H_i$. We define
$H_{i+1}$ by adding this $u$ and the outgoing edge.
Then $H_{i+1}$ clearly satisfies above two conditions.
Since $G$ is finite, we find $m$ that $G\setminus H_m$ is empty, i.e.,
the set of vertices of $G$ and $H_m$
are the same. We finish the proof by taking $G(C)=H_m$.
\qed
\bigskip

{\it Proof of Theorem \ref{Main}.}
Theorem 4.3
of \cite{Akiyama-Lee:13} shows that overlap coincidence of $\Tk$
implies multiple strong coincidence of level $n$ for any $n\ge 1$. We prove that
there is a constant $n$ such that multiple strong coincidence of level $n$ implies
overlap coincidence.

Assume that $\Tk$ does not admit overlap coincidence.
Construct the overlap graph $G$ of $\Tk$ with multiplicity.
Since $\Tk$ does not admit overlap coincidence, there is a strongly
connected component\footnote{In this assertion,
one can take either usual overlaps
or potential overlaps as we like.}
$S$ of $G$ such that its spectral radius
is equal to $|\det(Q)|$ and from each overlap of $S$
there is no path leading to a coincidence
in $G$.
Without loss of generality, we may assume that the incidence matrix of $S$ is primitive\footnote{If the incidence matrix of $S$ is irreducible but not primitive,
then take a suitable power of $\Omega$ by Perron-Frobenius theory.}.
Thus we can find a positive integer $n_0$
such that for every overlap $(U, y, V)$,
$\Omega^{n_0}(U, y, V)$ contains an overlap equivalent to
$(U, y, V)$.
Since $\Xi(\Tk)$ is a Meyer set, number of equivalence classes of
overlaps is finite and bounded by
a constant which depends only on $\Tk$.
Thus there is an upper bound of $n_0$ which depends only on $\Tk$.
We further assume multiple strong coincidence of level $n=n_0$ on $\Tk$
and derive a contradiction.

We claim
that in the component
$S$ there is an overlap $(U, y, V)$ with $\ell(U)\neq \ell(V)$
for any non-rod self-affine tiling by connected tiles.
Assume on the contrary that
all overlaps in $S$ are of the form $(A,z,B)$ with $\ell(A)=\ell(B)$.
Since $S$ does not contain a coincidence, $z\neq 0$ for these overlaps.
Taking $k$-th inflated overlap of $(A,z,B)$, we
obtain of patches $P$ and $Q$, both contain large balls, say $B_p(r)$ and $B_q(r)$,
that
the tiles of $P$ close to $p$ and the
tiles of $Q$ close to $q$ are in multiple correspondence in the following sense.
Putting $x=Q^k z$, for a tile $U\in P$ close to $p$
there are several (at least two) tiles $V\in Q$ that $(U,x,V)$ are overlaps in\footnote{We say that an overlap belongs to $S$ if its equivalence class does.}
 $S$
and $\supp(U)$ is contained in the union of $\supp(V-x)$,
and the same statements hold after interchanging the role of $U$ and $V$.
Take a tile $U$ with $p\in \supp(U)\subset B_p(r)$.
Then overlaps $(U,x,V)$ with $V\in Q$ give rise to a patch $\mathcal{V}_1
=\bigcup V$ that
$\supp(U) \subsetneq \supp(\mathcal{V}_1)-x$. By assumption, $\ell(V_1)=\ell(U)$
for every $V_1\in \mathcal{V}_1$.
By using path connectedness of tiles\footnote{Connectedness
and path connectedness are equivalent for self-affine tiles
\cite{Luo-Akiyama-Thuswaldner:04}.}, the patch $\mathcal{V}_1$ is path connected.
If $\supp(V_1)\subset B_q(r)$,
then there is a patch $\mathcal{U}_1=\bigcup U_1$ where $U_1\in P$ are taken
from all overlaps of the form $(U_1,x,V_1)$ with some $V_1\in \mathcal{V}_1$.
This patch
is also path connected and satisfies
$\supp(\mathcal{V}_1)-x \subsetneq \supp(\mathcal{U}_1)$
and each tile of $\mathcal{U}_1$ has the same color as $U$.
In this manner, by taking
large $r$, we obtain a long sequence of path connected patches
$$
\supp(U) \subsetneq \supp(\mathcal{V}_1-x) \subsetneq \supp(\mathcal{U}_1)
\subsetneq \supp(\mathcal{V}_2-x) \subsetneq \supp(\mathcal{U}_2) \subsetneq \dots.
$$
The number of tiles strictly increases and
all tiles appear in this sequence
has the same color $\ell(U)$.
This shows for any $M>0$, there exists a ball of radius $R$ that
each tile $U$ in the ball belongs to a connected patch in $\Tk$
having diameter greater than $M$, whose tiles have an identical color $\ell(U)$.
Therefore by using FLC, among $X_{\Tk}$ we can choose a rod tiling. Being
a rod tiling is invariant under
translation and closure operation,
using minimality of $X_{\Tk}$ we see that every tiling in $X_{\Tk}$
is a rod tiling. This gives a contradiction, which
finishes the proof of the claim.

Consider a directed
graph $\mathcal{V}$ over $\{1,\dots,m\}$
whose edge $i\rightarrow j$ is given if there are $U,V\in S$ that
$V\in \Omega^n(U)$ with $i=\ell(U)$ and $j=\ell(V)$. Clearly
$\mathcal{V}$ is strongly connected as well.
Pick one overlap $(U,y,V)$ from $S$ that $\ell(U)\neq \ell(V)$ and
select one of the overlaps equivalent to $(U,y,V)$ in $\Omega^n(U,y,V)$.
We select a tile map $\gamma=\gamma_{\Omega^{n}}$
which sends $\gamma(U)$ to this $U$ in $(U,y,V)$, and
$\gamma(V)$ to the $V$ in $(U,y,V)$, which correspond to
two cycles $\ell(U)\rightarrow \ell(U)$ and $\ell(V)\rightarrow \ell(V)$
on $\mathcal{V}$. Let $C$ be the set of these two cycles and take
$\mathcal{V}(C)$ by Lemma \ref{CycleExt}. The tile map
$\gamma=\gamma_{\Omega^n}$ is chosen so that $\ell(U)\rightarrow \ell(\gamma(U))$
for $U\in \{T_1,\dots,T_m\}$ forms the set of edges of $\mathcal{V}(C)$.
By the choice of the subgraph, every path of length $m$
on this subgraph must fall into one of the two cycles.
Note that by this choice of $\gamma$, the control points of $U$ and $V-y$ are
exactly matching, because both of them are equal to a
common point $\cap_{k=1}^{\infty} Q^{-nk}\left(\gamma^k(U)\cap \gamma^k(V-y)\right)$.

We claim that by this $\gamma$, we have $\Lam-\Lam \subset
\langle \mathcal{G} \rangle$.
In fact, since every overlap in the overlap graph is of the form $(A,z,B)$ with
$z\in \bigcup_{i=1}^m (\Lam_i-\Lam_i)$, and control points of $U$ and $V-y$ are matching on $(U,y,V)$,
i.e., $c(U)=c(V)-y$, we have
$c(U)-c(V) \in \mathcal{G}$.
By construction of $\mathcal{V}$
for any $x,y\in \Lam$, we have $Q^m x, Q^m y\in \Lam_{\ell(U)} \cup \Lam_{\ell(V)}$.
For instance, if $Q^m x\in \Lam_{\ell(U)}$ and
$Q^m y\in \Lam_{\ell(V)}$, then $Q^m x= c(U)+ f, Q^m y=c(V)+g$ hold
with $f\in \Lam_{\ell(U)}-\Lam_{\ell(U)},
g\in \Lam_{\ell(V)}-\Lam_{\ell(V)}$. Therefore we have $\Lam-\Lam\subset \langle
\mathcal{G} \rangle$.

We also see that the set of control points
$\Lam=(\Lam_i)$ associated to $\gamma$
is admissible. In fact, since $(U,y,V)$ is an overlap,  $\supp(U)\cap \supp(V-y)$
has an inner point. Since $y=c(V)-c(U)$, we have
$(\supp(U-c(U)))^{\circ}\cap (\supp(V-c(V)))^{\circ}\neq \emptyset$.
The admissibility follows from $Q^m x\in \Lam_{\ell(U)} \cup \Lam_{\ell(V)}$
for any $x\in \Lam$.

Summing up, from $\widetilde{(U,y,V)}\in S$,
we have chosen a tile map $\gamma_{\Omega^n}$ which produces
a substitution Delone multi-color set of
admissible control points with $\Lam-\Lam\subset \langle \mathcal{G} \rangle$.
By the assumption of multiple strong coincidence of level $n$, we know
$\Omega^k(U-c(U)) \cap \Omega^k(V-c(V))$
is non empty for some $k$, which shows that
$(U, y, V)$
leads to a coincidence, giving a desired contradiction.
\qed

\begin{remark}
\label{Rod}
We use the assumptions that each tile is connected and $\Tk$ is a non-rod tiling
only to show that there is an overlap $(U,y,V)\in S$ that $\ell(U)\neq \ell(V)$,
which allows us to define a tile map. It is likely that these assumptions are
not necessary, i.e., every non-periodic self-affine tiling that $\Xi(\Tk)$ is a
Meyer set, then such overlap must appear in $S$.
\end{remark}

\providecommand{\bysame}{\leavevmode\hbox to3em{\hrulefill}\thinspace}
\providecommand{\MR}{\relax\ifhmode\unskip\space\fi MR }
\providecommand{\MRhref}[2]{%
  \href{http://www.ams.org/mathscinet-getitem?mr=#1}{#2}
}
\providecommand{\href}[2]{#2}


\end{document}